\documentclass[reqno,11pt]{amsart}
\usepackage{pgfplots}

\usepackage{amsmath}
\usepackage{cancel}
\usepackage[normalem]{ulem}
\usepackage{breqn}
\usepackage{amssymb}
\usepackage{amsfonts}
\usepackage{mathtools}
\usepackage{enumerate}
\usepackage{hyperref}
\usepackage{comment}
\allowdisplaybreaks

\usepackage{palatino}
\usepackage[T1]{fontenc}
\usepackage[mathscr]{eucal}
\usepackage{dsfont}

\usepackage{fullpage}

\usepackage{acronym}
\usepackage{latexsym}
\usepackage{paralist}
\usepackage{wasysym}
\usepackage{xspace}

\numberwithin{equation}{section}  
\usepackage[sort&compress,capitalize,nameinlink]{cleveref}
\crefname{app}{Appendix}{Appendices}

\crefrangeformat{equation}{\upshape(#3#1#4)\textendash(#5#2#6)}

\theoremstyle{plain}
\newtheorem{theorem}{Theorem}
\newtheorem{corollary}[theorem]{Corollary}
\newtheorem*{corollary*}{Corollary}
\newtheorem{lemma}[theorem]{Lemma}
\newtheorem{proposition}[theorem]{Proposition}

\theoremstyle{definition}

\newtheorem*{definition*}{Definition}

\newtheorem*{hypothesis*}{Hypothesis}

\theoremstyle{remark}

\newtheorem*{remark*}{Remark}
\newtheorem*{notation*}{Notational remark}

\numberwithin{theorem}{section}

\DeclarePairedDelimiter{\abs}{\lvert}{\rvert}

\DeclarePairedDelimiter{\bra}{\langle}{\rvert}
\DeclarePairedDelimiter{\ket}{\lvert}{\rangle}
\DeclarePairedDelimiterX{\braket}[2]{\langle}{\rangle}{#1,#2}

\newcommand{\ind}{\mathds{1}}

\newcommand{\cI}{\ensuremath{\mathcal I}}

\newcommand{\cM}{\ensuremath{\mathcal M}}

\newcommand{\cP}{\ensuremath{\mathcal P}}

\newcommand{\cX}{\ensuremath{\mathcal X}}

\newcommand{\N}{\ensuremath{\mathbb{N}}}

\newcommand{\R}{\ensuremath{\mathbb{R}}}
\renewcommand{\P}{\ensuremath{\mathbb{P}}}

\newcommand{\pl}{\ensuremath{\left\langle}}
\newcommand{\pr}{\ensuremath{\right\rangle}}

\def\({\left(}
\def\){\right)}
\def\[{\left[}
\def\]{\right]}
\begin{document}
\title[]{Quasi-stationary distributions of non-absorbing Markov chains}

\author[R.~Fernandez]{Roberto Fernandez$^{\star}$}
\address{$^{\star}$ New York University Shanghai,  1555 Century Avenue, Pudong New Area, 200122, Shanghai, China}
\email{rf87@nyu.edu}

\author[F.~Manzo]{Francesco Manzo$^{*}$}
\address{$^{*}$ Dipartimento di Matematica e Fisica, Universit\`a di Roma Tre, Largo S. Leonardo Murialdo 1, 00146 Roma, Italy.}
\email{manzo.fra@gmail.com; scoppola@mat.uniroma3.it}

\author[M.~Quattropani]{Matteo Quattropani$^{\dagger}$}
\address{$^{\dagger}$ Dipartimento di Matematica, Universit\`a di Roma Tor Vergata, Via della Ricerca Scientifica 1, 00133, Roma, Italy.}
\email{quattropani@mat.uniroma2.it}

\author[E.~Scoppola]{Elisabetta Scoppola$^{*}$}

\maketitle              

\begin{abstract}
We consider reversible ergodic Markov chains with finite state space, and we introduce a new notion of quasi-stationary distribution that does not require the presence of any absorbing state. In our setting, the hitting time of the absorbing set is replaced by an optimal strong stationary time, representing the ``hitting time of the stationary distribution''. On the one hand, we show that our notion of quasi-stationary distribution corresponds to the natural generalization of the \emph{Yaglom limit}. On the other hand, similarly to the classical quasi-stationary distribution, we show that it can be written in terms of the eigenvectors of the underlying Markov kernel, and it is therefore amenable of a geometric interpretation. Moreover, we recover the usual exponential behavior that characterizes quasi-stationary distributions and metastable systems. We also provide some toy examples, which show that the phenomenology is richer compared to the absorbing case. Finally, we present some counterexamples, showing that the assumption on the reversibility cannot be weakened in general.
\end{abstract}
\tableofcontents

\section{Introduction}
In the theory of absorbing Markov chains, for a discrete-time chain $(X_t)_{t\ge 0}$, with state space $\cM\cup G$ where $G$ is the absorbing set, the \emph{quasi-stationary distributions} are the initial distributions supported on $\cM$ having the following special properties: on the one hand they exhibit an (exact) exponential absorption time and, on the other hand, they are invariant given that the chain has not been absorbed yet. These distributions can be characterized in two ways: in a linear-algebraic fashion, by looking at the principal eigenvectors of the Markov kernel restricted to $\cM$; or as the so-called \emph{Yaglom limits}, namely the limiting distributions, as $t\to\infty$, of the process conditioned to not have been absorbed at time $t$. 

Quasi-stationary distributions play a special role in metastability. 
In general, metastability has been associated with the existence of two timescales: one necessary to reach equilibrium and another, much shorter one, in which the system is well-described by a metastable measure. This metastable measure represents a kind of local equilibrium that is reached before stable equilibrium. As a striking feature of metastability, the exit from the metastable state occurs at an exponentially distributed time. 

\begin{figure}[h]
	\includegraphics[width=6cm]{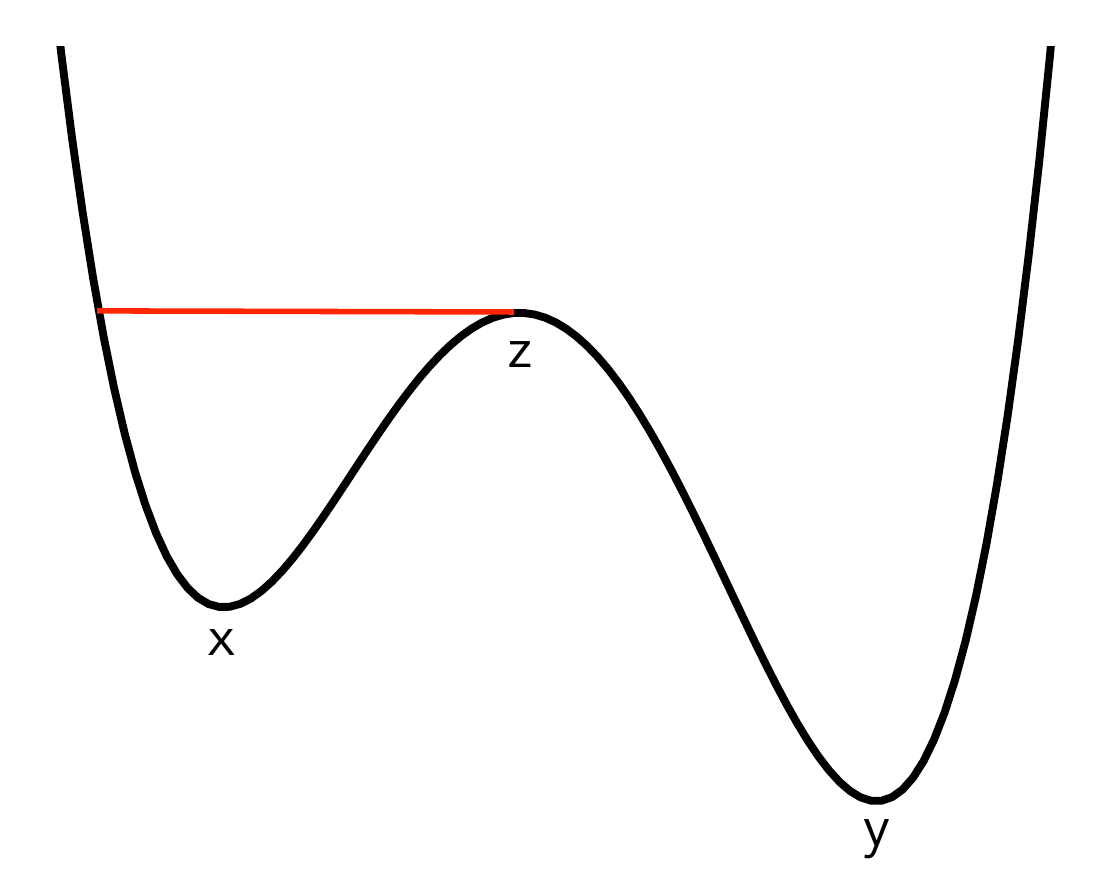}
	\caption{The picture shows the classical cartoon of a double-well energetic landscape in which the metastable state is associated with the bottom of the first well ($x$) whereas the stable state is associated to the bottom of the second well ($y$). The \emph{cycle} associated with the metastable state $x$ is the set of points lying below the red horizontal segment which is tangent the local maximum $z$ and intersects the left well.}\label{fig1}
\end{figure}

A classical setup in which metastability can be illustrated is the so-called Friedlin-Wentzell regime (Glauber-like dynamics on a finite set in the zero-temperature limit). In this regime the stationary and metastable measures are concentrated on  different regions. Metastability can then be described as the exit from special sets called \emph{cycles} which are defined in terms of the energy landscape and which include the metastable point. Finally, general results allow to show that as soon as the process leaves such metastable cycle, then the relaxation to the stable equilibrium occurs on a much shorter timescale. We refer the reader to \cref{fig1} for a classical cartoon illustrating the setup just described.

 	The study of metastability as an exit problem from a set establishes a connection between a metastable state and a quasi-stationary distribution, since the latter coincides with the long-time limit of the process conditioned not to reach the set's boundary, which is a rare event in Friedlin-Wentzell regime. As an immediate consequence of the identification of the metastable state with a quasi-stationary distribution, one has the exponential distribution of the exit time. In the early literature, where such a connection was not established yet, renewal theory was used for the same scope (see, e.g., \cite{OV05,MNOS,BdH}).

By looking at metastability as an exit problem from a connected set $\cM$,  the uniqueness of the Yaglom limit and its independence of the initial distribution $\alpha$ (as soon as it is supported only on $\cM$) immediately follow:
\begin{equation}\label{eq:vero-yaglom}
\lim_{t\to\infty} \P_\alpha\(X_t=y\mid \tau_G>t \)=\phi^*(y),\qquad y\in\cM,\quad  \alpha \in \cP(\cM)\,,
\end{equation}
where $\cP(\cM)$ denotes the set of probability distributions over $\cM$, $X_0\sim\alpha$, $\tau_G$ is the first hitting time to $G$ and $\phi^*$  is the principal left eigenvector of the Markov kernel restricted to $\cM$ 
\begin{equation}
\bra{\phi^*}[P]_ {\cM}=\lambda\bra{\phi^*}.
\end{equation}
Henceforth, starting at $\phi^*$, $\tau_G$ has a geometric distribution
\begin{equation}
\P_{\phi^*}(\tau_G>t)=\lambda^t.
\end{equation}

In \cite{MS19}, two of the authors used the language of strong stationary times to describe the exit problem, finding, in particular, an exact representation formula for the exit time distribution in terms of distributions of appropriate strong times.
The notion of Strong Stationary Time (SST) dates back to the works of Aldous and Diaconis \cite{AD86,AD2} (see \cite{LP} for a review): A SST is a randomized stopping time $\tau_{\pi}$ for which
$$  \P(X_t=x \mid \tau_{\pi}=t) = \pi(x)\,, $$
namely, the probability that the process $X$ is at $x$ at time $t$, conditioned on the event that SST is equal to $t$, is equal to the stationary measure $\pi$ evaluated at $x$. 
In this spirit, by introducing the notion of  \emph{conditional strong quasi-stationary time} the authors study  the convergence to the quasi-stationary distribution providing  an exact formula for the first hitting to the absorbing set $G$ with an explicit dependence on the starting distribution $\alpha$.

Outside of the Friedlin-Wentzell regime the analysis of metastable systems is much more subtle \cite{SS98,BGM20,BdH,DS97,CM}. Many results are still based  on studying a suitable exit problem, but this approach is useful only when the metastable state can be associated with a confined region of the state space, and in particular with the quasi-stationary distribution on this set.  
Therefore, the research of an analogue of the quasi-stationary measure to a non-absorbing setting could be seen as a first step to an effective approach to metastability in more general scenarios.

 In order to proceed toward such a generalization, we will need to find a proper counter-part of the hitting time of the absorbing set $G$. For each initial distribution $\alpha$, this role will be played by the so called \emph{optimal} SST which, in a sense, can be seen as the ``fastest'' SST associated to the starting distribution $\alpha$.
In what follows, we are going to show that the Yaglom limit associated to such a stopping time exists and, later, we will recover also a spectral interpretation of this non-absorbing quasi-stationary distribution in terms of the eigenvectors of the underlying ergodic kernel. Moreover, as we will see later, the dependence on the initial distribution $\alpha$ is more subtle compared to the absorbing case.  Indeed, by looking at the Fourier coefficients of an initial distribution, we can uniquely associate to it the corresponding Yaglom limit, partitioning the set of initial distributions into \emph{basins of attraction} of different Yaglom limits. \textcolor{black}{Some examples are discussed in \cref{sec:examples}, where also a 1-dimensional random walk is studied in detail.}

 Throughout the paper we stick to the assumption that the underlying Markov chain is reversible. This requirement is not merely of technical nature. Indeed, presenting some counterexamples, in the last section we show that our notion of Yaglom limit may fail to exists in a non-reversible setup, and the conditioned measure can show a periodic behavior in time, suggesting the presence of more general attractors. Therefore, the extension to non-reversible chains remains an open problem which we intend to investigate further in the future.

%
%
\section{Notation and results}
Consider a discrete-time Markov chain $(X_t)_{t\ge 0}$ on a finite state space $\cX$. In what follows, we will use the notation $\bra{\cdot}$ to refer to row vectors (measures) and $\ket{\cdot}$ to refer to column vectors (functions).
Assume that the associated transition matrix $P$ is irreducible and aperiodic, having the unique stationary distribution $\bra{\pi}$, which we assume to be reversible for $P$, and we further assume that the non-trivial eigenvalues of $P$ are all positive\footnote{This can be immediately obtained by considering the \emph{lazy} version of a general reversible chain $P$.}. Recall also that, by Perron-Frobenius theorem, all the entries of $\bra{\pi}$ are positive. It is worth to remark that all the results can be immediately extended to the reversible continuous-time setting. Nonetheless, in order to simplify the reading, we prefer to stick here to discrete-time.

Call $\bra{\alpha}$ the starting distribution of the chain and, for all $t\ge 0$, we use the notation $\bra{\mu_t^\alpha}$ to denote the distribution 
\begin{equation}
\mu_t^\alpha(y)\coloneqq 	\sum_{x\in\cX}\alpha(x)P^t(x,y),\qquad\forall y\in\cX\,.
\end{equation}
We will use the \emph{separation distance} to quantify the convergence to stationarity, i.e., we consider
\begin{equation}
	s^\alpha_t(y)\coloneqq 1- \frac{\mu^\alpha_t(y)}{\pi(y)},\qquad s^\alpha_t\coloneqq\max_{y\in\cX}s_t^\alpha(y)\,.
\end{equation}
Consider an \emph{optimal strong stationary time} $\tau_\pi^\alpha$, i.e., a randomized stopping time such that
\begin{equation}
	\P_\alpha(X_t=y,\: \tau_\pi^\alpha\le t)=\pi(y)\P_\alpha(\tau_\pi^\alpha\le t),\qquad \forall t\ge 0,\: y\in\cX\,,
\end{equation}
 furthermore, denote
\begin{equation}
\P_\alpha(\tau_\pi^\alpha> t)=s^\alpha_t,\qquad\forall t\ge0\,.
\end{equation}
Strong stationary time have been introduced in the pioneering work of Aldous and Diaconis \cite{AD86} and an explicit construction has been proposed by Diaconis and Fill in \cite{DF90}. For a modern introduction to the subject we refer to \cite[Ch. 6]{LP}.

We now consider the sequence of conditional distributions $(\bra{\varphi^\alpha_t})_{t\ge 0}$, defined as follows:
\begin{equation}\label{eq:cond-distr}
	\varphi^\alpha_t(y)\coloneqq\P_\alpha\(X_t=y\mid \tau_\pi^\alpha>t \),\qquad\forall t\ge 0,\:y\in\cX\,.
\end{equation}
The first result is a recursive formula for the sequence $(\varphi^\alpha_t)_{t\ge 0}$, which is contained in the following lemma.
\begin{lemma}\label{lemma:recursion-phi}
	For every $t\ge0$ it holds
	\begin{equation}\label{eq:recursion-phi}
		\bra{\varphi^\alpha_t}-\bra{\pi}=\frac{s_0^\alpha}{s_t^\alpha}\[\bra{\varphi^\alpha_0} P^t -\bra{\pi} \]\,,
				\end{equation}
	which in discrete-time rewrites as
		\begin{equation}\label{iter}
			\bra{\varphi^\alpha_t}=\frac{s_{t-1}^\alpha}{s_t^\alpha}\bra{\varphi^\alpha_{t-1}}P +\(1-\frac{s_{t-1}^\alpha}{s_t^\alpha}\)\bra{\pi}\,.
	\end{equation}
\end{lemma}
As a corollary of \cref{lemma:recursion-phi} we show that the sequence of distributions $(\varphi^\alpha_t)_{t\ge0}$ are \emph{maximally} separated from $\pi$, in the sense that they cannot be fully supported on $\cX$. In what follows, we will call \emph{simplex} the (compact) set of probability distributions on $\cX$.
\begin{corollary}\label{coro:sepa}
For every $t\ge0$ it holds
\begin{equation}
	\sup_{y\in \cX} 1-\frac{\varphi^\alpha_t(y)}{\pi(y)}=1\,,
\end{equation}
 i.e., for every $t\ge0$ there exists some $y\in\cX$ such that $\varphi^\alpha_t(y)=0$.
\end{corollary}
In words, \cref{lemma:recursion-phi} and \cref{coro:sepa} tell us that the sequence of conditional measures $(\varphi_{t}^\alpha)_{t\ge 0}$ evolves by remaining on the border of the simplex and, at each step, the conditional measure is a convex combination of the evolved conditional measure at the previous step  and the stationary distribution. To interpret such phenomenology we notice that,
in the spirit of \cite{FL14}, the optimal strong stationary time of the process $P$ started at $\alpha$ has the same law of a hitting time of a countable-state Markov chain. Indeed, consider the Markov chain on $\N_0\cup\{\dagger \}$ with transition matrix
	\begin{equation}
		\hat P_\alpha(t,t')=\begin{cases}
			1-\frac{s^\alpha_{t+1}}{s^\alpha_t}&t\neq\dagger\,,t'=\dagger\,,\\
			\frac{s^\alpha_{t+1}}{s^\alpha_t}&t\neq\dagger\,,t'=t+1\,,\\
			1&t=\dagger\,,t'=\dagger\,,\\
			0&\text{otherwise}\,,
		\end{cases}
	\end{equation}
	and starting at $0$.  In other words, $\hat P_\alpha$ describes a pure-birth process with killing: at each step the chain can either do a step forward on the positive integer line, or get killed. 
	Thanks to \cref{lemma:recursion-phi}, one gets the intertwining relation
	\begin{equation}
		\Lambda_\alpha P=\hat P_\alpha\Lambda_\alpha
	\end{equation}
 where, for every $\alpha\neq\pi$,
	\begin{equation}
		\Lambda_\alpha(t,x)=
		\begin{cases}
			\varphi^\alpha_t(x)& t\ge 0\,,x\in\cX\,,\\
			\pi(x)&t=\dagger\,,x\in\cX\,.
		\end{cases}
	\end{equation}
Choosing as initial distribution for the chain $\hat P$ the distribution $\hat\alpha\in\cP(\N_0\cup\{\dagger\})$ defined as
	\begin{equation}
		\hat\alpha=s^\alpha_0 \delta_0 + (1-s^\alpha_0)\delta_\dagger\,,
	\end{equation}
	we have
	\begin{equation}
		\P_\alpha(\tau^\alpha_{\pi}=t)=\hat{\P}_{\hat{\alpha}}(\tau_\dagger=t)\,,\qquad t\ge 0\,.
	\end{equation}
Hence, one can interpret the sequence $(\varphi_{t}^\alpha)_{t\ge 0}$ as the sequence of ``states'' visited by the chain up to the occurrence of the SST, when the process eventually reaches its equilibrium.

Under this perspective, one could say that metastability occurs if there exists a timescale $t$ for which  $\P(\tau_{\pi}^\alpha>t)\approx 1$ and $\varphi^\alpha_s\approx \varphi^\alpha_{s+1}$ for all $s\ge t$.  Therefore, as a first step, it is natural to ask under which circumstances the limiting conditional distribution $\lim_{t\to\infty}\varphi_{t}^\alpha$ exists and, in the affirmative case, to interpret it as a candidate \emph{metastable state}.
The next theorem show that reversibility suffices to guarantee the existence of the latter limit, simultaneously for all $\alpha\neq \pi$.
\begin{theorem}[Existence of the Yaglom limit]\label{conj:M0}
	If the transition matrix $P$ is reversible, irreducible and with positive eigenvalues, then for every starting distribution $\alpha\neq\pi$ there exists the limit
	\begin{equation}\label{eq:def-limit}
		\varphi^\alpha_\star(y)\coloneqq\lim_{t\to\infty} \P_\alpha\(X_t=y\mid \tau^\alpha_\pi>t \),\qquad y\in\cX.
	\end{equation}
\end{theorem}
The proof of \cref{conj:M0} exploits the reversibility of the transition matrix $P$,  revolving around its spectral decomposition. Our proof of \cref{conj:M0} is constructive, in the sense that we derive an operative formula for the limit in \cref{eq:def-limit} in terms of $\alpha$ and the eigenvectors and eigenvalues of $P$. Such an explicit formula for $\varphi_{\star}^\alpha$ can be found in \cref{eq:varphistar_general}. Since we believe that the formula itself, in its full generality, does not favor much the intuition, in the next statement we provide a formula which holds true only under stronger assumptions, but which is easier to read.
\begin{corollary}\label{coro:simplified}
Let $P$ as in \cref{conj:M0} and assume further that all the eigenvalues of $P$ are distinct. Let $1=\lambda_1>\lambda_2>\dots>\lambda_n> 0$ the eigenvalues of $P$ and $\bra{{v}_1},\dots,\bra{{v}_n}$ the associated left eigenvectors. For every $\alpha\neq\pi$, call $\lambda^\alpha$ the largest eigenvalue such that the associated left eigenvector satisfies $\braket{{v}^\alpha}{\alpha}\neq 0$. Then
\begin{equation}\label{eq:varphistar_simplified}
	\bra{\varphi_{\star}^\alpha}= \bra{\pi}+\left(\frac{1}{\max_{y\in\cX_-^\alpha}\left|v^\alpha(y)\pi(y)^{-1}\right|}\right)\bra{{v}^\alpha} \,,
\end{equation}
where 
\begin{equation}\label{cX}
	\cX_-^\alpha\coloneqq\{y\in\cX\mid v^\alpha(y)<0 \}\,.
\end{equation}
\end{corollary}
The general formula is not much different than that in \cref{coro:simplified}. Nevertheless, rather than write it explicitly here, we believe it is more convenient to collect the main features of the limiting distribution in \cref{eq:def-limit} in the next result, making an explicit connection between the Yaglom limit and the spectral properties of the associated transition matrix.
\begin{corollary}\label{coro:consequences}
	If the limiting distribution in \cref{conj:M0} exists, then for every $\alpha\neq\pi$ there exists a couple $(v^\alpha,\lambda^\alpha)\in \R^{|\cX|}\times (0,1)$ such that  $(v^\alpha,\lambda^\alpha)$ satisfy
	\begin{equation}\label{eq:evector-v}
		\bra{v^\alpha} P = \lambda^\alpha \bra{v^\alpha}\,,
	\end{equation}
and
	\begin{enumerate}
		\item The distribution $\varphi_\star^\alpha$ is in the span of $v^\alpha$ and $\pi$. In particular
		\begin{equation}\label{eq:evector-phi}
			\bra{\varphi_\star^\alpha} = \bra{v^\alpha}+\bra{\pi}\,.
		\end{equation}
		\item The distribution $\varphi_\star^\alpha$  satisfies
		\begin{equation}\label{eq:evector}
			\bra{\varphi_\star^\alpha} P^t= (\lambda^\alpha)^t	\bra{\varphi_\star^\alpha} + \big(1- (\lambda^\alpha)^t \big) \bra{\pi}\,, \qquad \forall t\ge0\,.
		\end{equation}
		\item The separation distance, starting at $\alpha$, decays exponentially at rate $\lambda^\alpha$, i.e.,
		\begin{equation}\label{eq:decay-separation}
			\lim_{t\to\infty}\(s^\alpha_t \)^{\frac{1}{t}}=\lambda^\alpha\,.
		\end{equation}
	\item It exists some state $y\in\cX$ such that
	\begin{equation}\label{eq:halting-alpha}
	\varphi_\star^\alpha(y)=0\,.
	\end{equation}
	\end{enumerate}
\end{corollary}
In words, \cref{coro:consequences} shows that our notion of Yaglom limit exhibits all of the properties one would expect from a metastable measure: it can be represented in terms of the eigenvectors of $P$ (property (1)); starting at the Yaglom limit, the (unconditional) process remains in a mixture of the stationary distribution and the Yaglom limit itself, and the decay to equilibrium is exponential (property (2)). Moreover, property (3) provides an explicit decay rate for the separation distance to equilibrium starting at $\alpha$, while property (4) tells us that the Yaglom limit is always maximally separated from the equilibrium.
\subsection{Geometric interpretation}
\begin{figure}[h]
	\includegraphics[width=10cm]{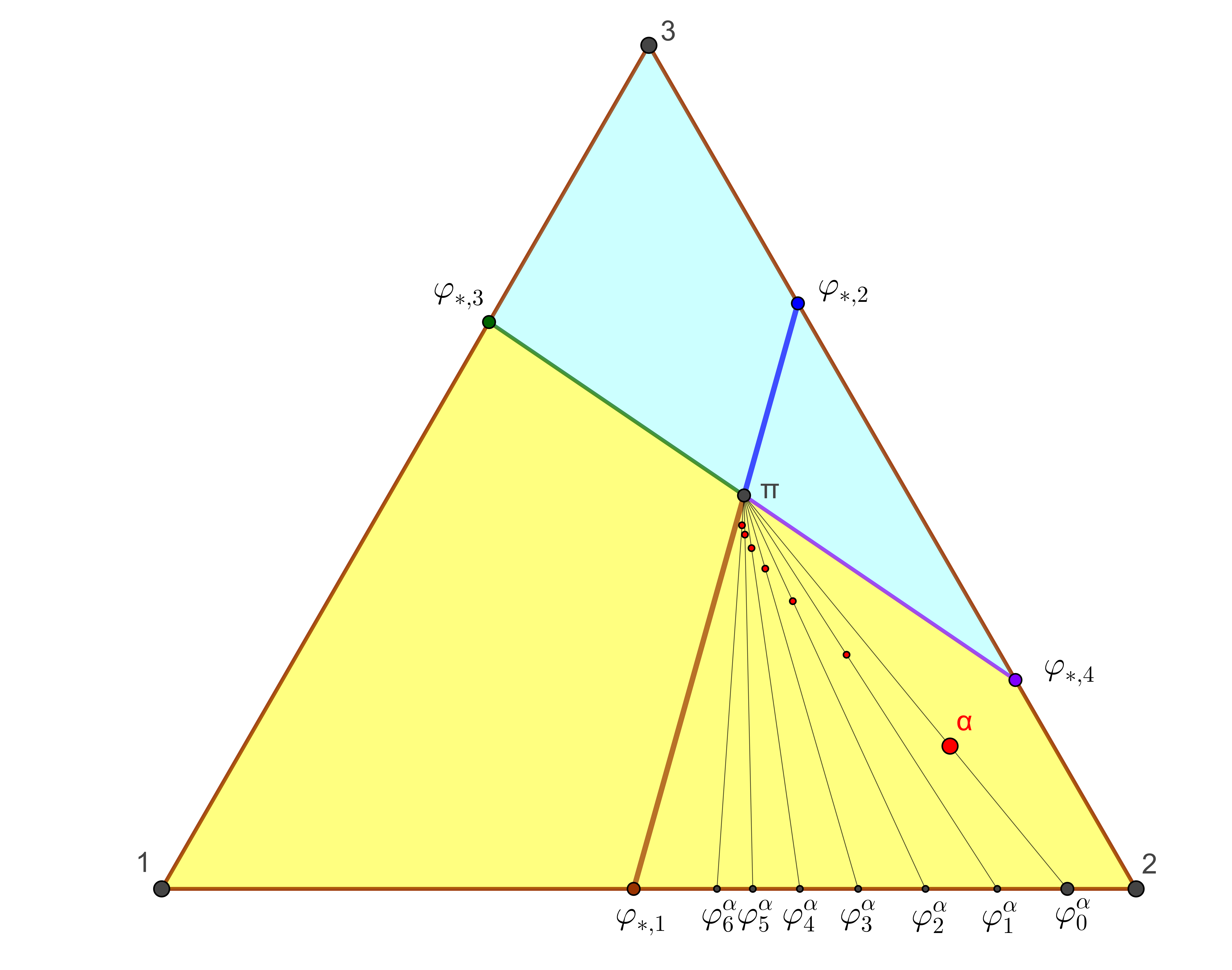}
	\caption{Geometric representation of the 2-dimensional simplex together with the evolution of the unconditional, $(\mu^\alpha_t)_{t\ge 0}$, and conditional,  $(\varphi^\alpha_t)_{t\ge 0}$ distributions, in red and in black, respectively. In particular, the vector $v_2$ (resp. $v_3$) is parallel to the segment joining $\varphi_{\star,1}$ and $\varphi_{\star,2}$ (resp. $\varphi_{\star,3}$ and $\varphi_{\star,4}$ ). This picture is a numerical simulation of the concrete toy example analyzed in detail in \cref{sec:1drw}.}\label{fig:cartoon}
\end{figure}

The notion of quasi-stationary distribution introduced in this work allows  to describe the way in which a stochastic process approaches its equilibrium state through the probabilistic lens of strong stationary times. Such a description has a linear-algebraic counterpart, which is naturally amenable of a geometric interpretation. Let us illustrate further this interpretation by exploiting a simple 3-state example.
In \cref{fig:cartoon}, we represent the 2-dimensional simplex, i.e., the compact set of probability distribution over 3-states. The vertices of the simplex represent the Dirac distributions on the three states, where the edges coincide with the set of probability distributions which are supported on only two states (those associated to the vertices at their extremes). \textcolor{black}{Therefore, for a given transition matrix $P$ for which one of the three states is absorbing, the stationary distribution $\pi$
coincides with one of the vertices of the simplex. In our non-absorbing (and irreducible) setup, instead, $\pi$ is represented by a point in the interior of the simplex.} By this perspective, is immediate that any convex combination of two distributions can be represented as a point lying on the segment joining them. Therefore, given the transition matrix $P$ (and thus $\pi$) and any distribution $\mu$, it is always possible  to rewrite $\mu$ as a convex combination of $\pi$ a third distribution $\phi$, i.e.,
\begin{equation}\label{eq:convex}
	\mu=(1-s)\pi+s\phi\,,
\end{equation}
where $s\in[0,1]$ and $\phi$ is another point on the simplex. Of course, there are infinitely many choices for the couple $(s,\phi)$ for which \cref{eq:convex} holds true. A natural way to choose one couple satisfying \cref{eq:convex} is by minimizing the value of $s$. We will refer to this choice with the name of \emph{greedy choice}, and we let $(s_\star,\phi_\star)$ denote the couple. Inverting the relation in \cref{eq:convex} we get
$$s=\frac{\pi(x)-\mu(x)}{\pi(x)-\phi(x)}\,,\qquad \forall x\in\cX\,.$$
In words, given $\phi$, $s$ represents the ratio between the length of the segment joining $\pi$ and $\mu$, and the length of the segment joining $\pi$ and $\phi$. Therefore, it is immediate to check that $\phi_\star$ must be a point on the border of the simplex and that $s_\star$ coincides with the separation distance between $\mu$ and $\pi$. In conclusion, for any $\mu$ the associated \emph{greedy} choice for $\phi_\star$ is nothing more than the unique point on the border of the simplex intersecting the half-line starting at $\pi$ and passing through $\mu$, while $s={\rm sep}(\mu,\pi)$. In this language, we deduce that the sequence of conditional distribution $\varphi_{t}^\alpha$ defined in \cref{eq:cond-distr} is noting else that the associated sequence of \emph{greedy choices} associated to $\mu^{\alpha}_t$. Therefore, it is clear that $(\varphi_{t}^\alpha)_{t\ge 0}$ describes a discrete-time dynamical system living on the border of the simplex. 

Let us stick again to the 3-state example represented in \cref{fig:cartoon}. Assuming that the process is reversible, then the two (non-principal) left eigenvectors can be represented as two lines on the plane containing the simplex, which cross each other at $\pi$. In order to aid the intuition, let us assume that the two associated eigenvalues do not coincide, i.e., $\lambda_2<\lambda_3$, and let us call $v_2$ and $v_3$ the associated lines. Clearly $v_3$ divides the simplex into two regions, and each $\alpha$ lying exactly on the line associated to $v_3$ is orthogonal to the right eigenvector associated to $v_2$, say $f_2$. Therefore, if $\alpha$ lies on $v_3$, the associated Yaglom limit (and actually the whole sequence of conditional distributions in \cref{eq:cond-distr}) coincide with the intersection between the border of the simplex and the half-line starting at $\pi$ and passing through $\alpha$. More interestingly, if $\alpha$ is not a point on $v_3$, then the associated Yaglom limit coincides with the unique intersection point between $v_2$ and the border of the simplex which is on the same side as $\alpha$ with respect to $v_3$. Hence, we have in total four possible Yaglom limits, each having its own \emph{basin of attraction}. Nevertheless, only two basins (those associated with Yaglom limits that are on the span of $\pi$ and $v_2$) have the same dimension of the simplex, while the other two have a lower dimensionality.

\section{Proofs}
\begin{proof}[Proof of \cref{lemma:recursion-phi}]
	Let us start by rewriting, for some $y\in\cX$ fixed,
	\begin{equation}\label{eq:rec}
		\begin{split}
			\varphi^\alpha_t(y)-\pi(y)&=\frac{\P_\alpha(X_t=y)-\P_\alpha(X_t=y,\: \tau_\pi^\alpha\le t)}{\P_\alpha(\tau_\pi^\alpha>t)}-\pi(y)\\
			&=\frac{\P_\alpha(X_t=y)-\pi(y)(1-s_t^\alpha)-s_t^\alpha \pi(y)}{s_t^\alpha}\\
			&=\frac{s_{0}^\alpha}{s_t^\alpha}\(\frac{\P_\alpha(X_t=y)-\pi(y)}{s_{0}^\alpha}\).
		\end{split}
	\end{equation}
	We now notice that
	\begin{equation}\label{eq:alpha}
		\alpha(y)=s_0^\alpha \varphi^\alpha_{0}(y) + (1-s_0^\alpha)\pi(y)\,.
	\end{equation}
	Applying $P^t$ to both sides of \cref{eq:alpha} and rearranging the terms we get
	\begin{equation}\label{eq:phi0Pt}
		\varphi^\alpha_{0}P^t(y)-\pi(y)=\frac{\P_\alpha(X_t=y)-\pi(y)}{s_0^\alpha}\,.
	\end{equation}
	Plugging \cref{eq:phi0Pt} into \cref{eq:rec} we deduce the desired claim.
\end{proof}
\begin{proof}[Proof of \cref{coro:sepa}]
	Rearranging the terms in \cref{eq:recursion-phi} we obtain
	\begin{equation}\label{eq:rearranging}
		s_t^\alpha \(\bra{\varphi^\alpha_t}-\bra{\pi} \)=s_0^\alpha \(\bra{\varphi^\alpha_{0}} P^t-\bra{\pi} \)\,.
	\end{equation}
	Looking at \cref{eq:rearranging} component-wise, and dividing by $\pi(\cdot)$ we obtain
	\begin{equation}\label{eq:rearranging2}
		\(1-\frac{\varphi^\alpha_t(y)}{\pi(y)} \)=\frac{s_0^\alpha}{s_t^\alpha} \(1 -\frac{\varphi^\alpha_{0} P^t(y)}{\pi(y)} \)\,.
	\end{equation}
	By using again \cref{eq:phi0Pt} it follows that
	\begin{equation}\label{eq:phi0Pt2}
		\sup_{y\in\cX}1-\frac{\varphi^\alpha_{0}P^t(y)}{\pi(y)}
		=\frac{1}{s_0^\alpha}	\sup_{y\in\cX}\[  1-\frac{\mu_t^\alpha(y)}{\pi(y)} \] =\frac{s_t^\alpha}{s_0^\alpha}\,.
	\end{equation}
	The claim follows by taking the supremum over $y\in\cX$ in \cref{eq:rearranging2} and plugging \cref{eq:phi0Pt2} on the right-hand side.
\end{proof}

\begin{proof}[Proof of \cref{conj:M0} and \cref{coro:consequences,coro:simplified}]
	Notice that, plugging  \cref{eq:phi0Pt} into \cref{eq:recursion-phi}, we obtain
	\begin{equation}\label{eq:final-expression}
		\bra{\varphi_t^\alpha}=\frac{1}{s_t^\alpha}\left(\bra{\mu_t^\alpha}-\bra{\pi}\right)+\bra{\pi}\,.
	\end{equation}
Therefore, the existence of the limit is equivalent to the existence of the limiting vector
\begin{equation}\label{eq:limiting-ratio}
\lim_{t\to\infty} \frac{1}{s_t^\alpha}\left(\bra{\mu_t^\alpha}-\bra{\pi}\right)\,.
\end{equation}
Call $n=|\cX|$ and rewrite
\begin{equation}
	P=\ket{1}\bra{\pi}+\sum_{i=2}^n\lambda_i \ket{f_i} \bra{v_i},\qquad \bra{\alpha}=a_1\bra{\pi}+\sum_{i=2}^n a_i \bra{v_i}\,,
\end{equation}
\begin{color}{black}
where $(\bra{v_i})_{2\le i\le n}$ denotes the collection of left eigenvectors of $P-\ket{1}\bra{\pi}$ and $(\ket{f_i})_{2\le i\le n}$ the corresponding right eigenvectors, normalized so to have
\begin{equation}
	\braket{v_i}{f_j}=\begin{cases}
		1&\text{ if }i=j\\
		0&\text{ if }i\neq j\,.
	\end{cases}
\end{equation}
In words, $a_i=\pl\alpha,f_i\pr$ is the $i$-th coordinate of $\alpha$ in the basis of the left eigenvectors.
\end{color}
Notice that, since $\alpha$ is a probability distribution, it must hold
\begin{equation}
	1=\braket{\alpha}{1}=\sum_{i=1}^n a_i \braket{v_i}{1}=a_1\,,
\end{equation}
therefore
\begin{equation}\label{eq:alpha-decomp}
	\alpha=\bra{\pi}+\sum_{i=2}^n a_i \bra{v_i}\,.
\end{equation}
It follows that
\begin{equation}\label{eq:mu-t-alpha}
	\bra{\mu_t^\alpha}-\bra{\pi}= \sum_{i=2}^n \lambda_i^ta_i \bra{v_i}\,.
\end{equation}
Recall that by assumption we have $1=\lambda_1>\lambda_2\ge \dots\ge \lambda_n> 0$, call
\begin{equation}
\lambda^\alpha\coloneqq \max_{i\in\{2,\dots,n \}}\{\lambda_i\text{ s.t. }a_i\neq 0 \}\,,
\end{equation}
and let
\begin{equation}\label{eq:def-cI-alpha}
\cI^\alpha\coloneqq \{i\in\{2,\dots,n \}\text{ s.t. }\lambda_i=\lambda^\alpha \}\,.
\end{equation}
Then, by \cref{eq:mu-t-alpha} and the definition of $\lambda^\alpha$ we deduce
\begin{equation}
	\lim_{t\to\infty} \frac{\bra{\mu_t^\alpha}-\bra{\pi}}{(\lambda^\alpha)^t}=\sum_{i\in \cI^\alpha}a_i\bra{v_i}\,.
\end{equation}
Therefore, in order to prove \cref{conj:M0} it is enough to show the existence of the limit
\begin{equation}\label{eq:limit}
	\ell^\alpha\coloneqq\lim_{t\to\infty} \frac{s_t^\alpha}{(\lambda^\alpha)^t}\,.
\end{equation}
If the limit exists, then
\begin{equation}\label{eq:phi-star-ell}
	\bra{\varphi^\alpha_\star}=\frac{1}{\ell^\alpha}\sum_{i\in \cI^\alpha}a_i\bra{v_i}+\bra{\pi}\,.
\end{equation}
We now show that the limit in \cref{eq:limit} exists. By \cref{eq:mu-t-alpha} we have
\begin{equation}
	\mu_t^\alpha(y)= \sum_{i=2}^n \lambda_i^ta_i v_i(y) + \pi(y)\,,
\end{equation}
hence,
\begin{equation}
	1-\frac{\mu_t^\alpha(y)}{\pi(y)}=- \sum_{i=2}^n \lambda_i^ta_i \frac{v_i(y)}{\pi(y)}\,.
\end{equation}
It follows that
\begin{equation}
	s_t^\alpha=\max_{y\in\cX}\sum_{i=2}^n -\lambda_i^ta_i \frac{v_i(y)}{\pi(y)}\,.
\end{equation}
Dividing both sides by $(\lambda^\alpha)^t$ we finally get
\begin{equation}
	\frac{s_t^\alpha}{(\lambda^\alpha)^t}=\max_{y\in\cX}\sum_{i=2}^n -\(\frac{\lambda_i}{\lambda^\alpha}\)^ta_i \frac{v_i(y)}{\pi(y)}\,.
\end{equation}
Notice that, when $t\to\infty$, we are left with
 \begin{equation}\label{eq:limit-sep}
 	\lim_{t\to\infty}\frac{s_t^\alpha}{(\lambda^\alpha)^t}=\max_{y\in\cX}\sum_{i\in \cI^\alpha} -a_i \frac{v_i(y)}{\pi(y)}\,,
 \end{equation}
so that
\begin{equation}\label{eq:ell-alpha}
	\ell^\alpha= \max_{y\in\cX}\sum_{i\in \cI^\alpha} -a_i \frac{v_i(y)}{\pi(y)}\,,
\end{equation}
and, by \cref{eq:phi-star-ell}, we deduce
\begin{equation}\label{eq:varphistar_general}
		\bra{\varphi^\alpha_\star}=\frac{1}{\max_{y\in\cX}\sum_{i\in \cI^\alpha} -a_i \frac{v_i(y)}{\pi(y)}}\sum_{i\in \cI^\alpha}a_i\bra{v_i}+\bra{\pi}\,.
\end{equation}
Notice that, in the case in which $|\cI^\alpha|=1$ and the coordinate along that axis is $a$, then the latter expression boils down to the claim of \cref{coro:simplified}.

To prove \cref{coro:consequences}, it is enough to realize that properties (1)-(2) are immediate consequences of \cref{eq:varphistar_general} and the definition of $\cI^\alpha$ in \cref{eq:def-cI-alpha} and property (3) follows from \eqref{eq:limit-sep}. We are left to show the validity of property (4). To this aim, call 
\begin{equation}
	\bra{v}=\sum_{i\in\cI^\alpha}a_i\bra{v_i}\,,
\end{equation}
and
	\begin{equation}\label{eq:vbar}
	\bar v(y) = \frac{v(y)}{\max_{z\in\cX_-^\alpha}\abs{\frac{v(z)}{\pi(z)}}}\,,\qquad  y\in\cX\,,
\end{equation} 
where $	\cX^\alpha_-$ is defined as in \cref{cX}, and notice that $\bar v$ is exactly the first vector on the right-hand-side of \cref{eq:varphistar_general} since
\begin{equation}
	\max_{y\in\cX}- \frac{v(y)}{\pi(y)}=\max_{y\in\cX_v^-} \left|\frac{v(y)}{\pi(y)} \right|\,.
\end{equation}
Hence, 
\begin{equation}
	\bra{\varphi_\star^\alpha}= \bra{\bar v}+\bra{\pi}\,.
\end{equation}
\textcolor{black}{To prove that $\varphi_{\star}^\alpha$ is indeed a probability distribution}, notice that by orthogonality
\begin{equation}
	\sum_{y\in\cX} \varphi_\star^\alpha(y)=1\,.
\end{equation}
and
\begin{equation}
	\varphi_\star^\alpha(y)\ge 0\,,\qquad y\in\cX\,.
\end{equation}
To see the latter we argue by contradiction, assuming that there exists some $y\in\cX$ such that $\varphi_\star^\alpha(y)<0$. Then
\begin{equation}
	0> \varphi_\star^\alpha(y)= \frac{v(y)}{\max_{z\in\cX^\alpha_-}\abs{\frac{v(z)}{\pi(z)}}} + \pi(y)\,,
\end{equation}
which leads to the contradiction
\begin{equation}
	-\frac{\pi(y)}{v(y)}> \max_{z\in\cX^\alpha_-} \left|\frac{\pi(z)}{v(z)}\right|\,.
\end{equation}
At this point, to conclude the proof of property (4) it is enough to call
\begin{equation}
	\tilde y\coloneqq \arg\max_{z\in\cX^-_v} \left|\frac{v(z)}{\pi(z)}\right|\,,
\end{equation}
and notice that
\begin{equation}
	\varphi_\star^\alpha(\tilde y)= \frac{v(\tilde y)}{\abs{\frac{v(\tilde y)}{\pi(\tilde y)}}} + \pi(\tilde y)=0\,.
\end{equation}
\end{proof}

\section{Examples}\label{sec:examples}
In this section we present some toy examples  in which the Yaglom limits can be explicitly computed. We start in \cref{suse:ex1} by studying what is arguably the easiest non trivial example: an homogeneous three-state chain. Even in this simple scenario, we show that the form of the Yaglom limit is not completely trivial. We then focus on two examples which present the features of a metastable system: in \cref{suse:ex2} we look at a two-block graph with a bottleneck, i.e., a chain presenting a metastability of entropic nature; on the other hand, in \cref{sec:1drw}, we study a three-state 1-dimensional random walk in the Friedlin-Wentzell regime (i.e., the most classical example of energetic metastability), and compare our notion of Yaglom limit with the classical (absorbing) one.
\subsection{The homogeneous triangle}\label{suse:ex1}
We start by considering the first non trivial example of a mean-field chain with $3$ states.
More precisely, we consider the transition matrix
\begin{equation}
	P=\frac{1}{4}\left(\begin{matrix}
		2&1&1\\
		1&2&1\\
		1&1&2
	\end{matrix} \right)\,,
\end{equation} 
having uniform stationary distribution $\pi\equiv\frac{1}{3}$. The eigenvalues are $(1,1/4,1/4)$ and a system of orthonormal eigenvectors spanning the vector space orthogonal to $\pi$ is given by
\begin{equation}
	\nu_2=\frac{1}{2\sqrt{3}}\(\sqrt{3}-1,2,-\sqrt{3}-1 \),\qquad 	\nu_3=\frac{1}{2\sqrt{3}}\(-\sqrt{3}-1,2,\sqrt{3}-1 \)\,.
\end{equation}
Consider an initial distribution $\bra{\alpha}=(\alpha(1),\alpha(2),\alpha(3))$. Writing $\bra{\alpha}$ in \cref{eq:alpha-decomp} we obtain
\begin{equation}\label{eq:alpha-triangle}
	\bra{\alpha}=\bra{\pi}+a_2\bra{\nu_2}+a_3\bra{\nu_3}\,,
\end{equation}
where
\begin{equation}\label{eq:c-triangle}
	a_2=\frac{1}{2}\alpha(1)+\frac{\sqrt{3}}{2}\alpha(2)-\frac{1}{2}\alpha(3)-\frac{1}{2\sqrt{3}}\,,\qquad
	a_3=-\frac{1}{2}\alpha(1)+\frac{\sqrt{3}}{2}\alpha(2)+\frac{1}{2}\alpha(3)-\frac{1}{2\sqrt{3}}\,.
\end{equation}
Assuming $\alpha\neq\pi$ we necessarily have $\lambda^\alpha=\frac{1}{4}$ and an immediate algebraic manipulation shows that, called $\alpha_{\min}\coloneqq\min\{\alpha(1),\alpha(2),\alpha(3) \}$,
\begin{equation}\label{eq:ell-triangle}
	\ell^\alpha=1-3\alpha_{\min}\,.
\end{equation}
Plugging \cref{eq:c-triangle,eq:ell-triangle} into \cref{eq:phi-star-ell}, we conclude that
\begin{equation}
\bra{\varphi_\star^\alpha}=\frac{1}{1-3\alpha_{\min}}\bigg(\alpha(1)-\alpha_{\min},\alpha(2)-\alpha_{\min},\alpha(3)-\alpha_{\min}\bigg)\,.
\end{equation}
In other words, the Yaglom limit $\varphi_\star^\alpha$ is obtained, starting from $\alpha$ and reducing each entry as soon as one (or more) of them becomes zero. The non-negative vector obtained in this way is then normalized in order to make it a probability vector. Notice that in this framework the set of Yaglom limits coincides with the boundary of the simplex (seen as a $2$-dimensional manifold in $\R^3$). Moreover, the basin of attraction of each Yaglom limit is a $1$-dimensional segment on the simplex.
	\subsection{An inhomogeneous mean-field example}\label{suse:ex2}
	Consider the case of a Markov chain with $2n$ states divided in two classes $\cX=\cX_1\sqcup\cX_1$ of the same size $n$. Assume that the transition matrix is of the form
	$$Q(x,y)=\begin{cases}
		\frac{1-\varepsilon}{n}&\text{if }x,y\in \cX_i\text{ for }i\in\{1,2\},\\
		\frac{\varepsilon}{n}&\text{otherwise},
	\end{cases} $$
	for some $\varepsilon\in(0,1/2)$. In other words, called $J$ the $n\times n$ matrix having $1/n$ in all entries, $A=(1-\varepsilon)J$ and $B=\varepsilon J$ then
	$$Q=\left(\begin{matrix}
		A&B\\
		B&A
	\end{matrix} \right).$$
	To avoid the existence of zero-eigenvalues consider, for some $\gamma\in(0,1)$ the lazy version
	$$P=(1-\gamma) Q+\gamma I.$$
	The eigenvalues of such a transition matrix are
	$\lambda_1=1$, $\lambda_2=(1-\gamma)(1-2\varepsilon)+\gamma$ and $\lambda_3=\gamma$ where the latter has multiplicity $2n-2$. Moreover, the eigenvector associated to $\lambda_2$ is such that $v(x)=v(y)$ for each $x,y\in\cX_i$ for some $i\in{1,2}$ and $v(x)=-v(y)$ if there is no $i\in\{1,2\}$ such that $x,y\in\cX_i$. Let us make the normalized choice $v(x)=\frac{1}{\sqrt{2n}}\ind_{x\in\cX_1}-\frac{1}{\sqrt{2n}}\ind_{x\in\cX_2}$. It is also immediate to check that the $2n-2$ vectors
	\begin{align*}
		w_i(x)=&\ind_{x=1}-\ind_{x=i+1}\,,\quad z_i(x)=\ind_{x=n+1}-\ind_{x=n+i+1}\,,\qquad i\in\{1,\dots n-1 \}\,,
	\end{align*} 
	are eigenvectors associated to $\lambda_3=\gamma$. We point out that the  $w$'s vectors (similarly for the $z$'s), even if not orthogonal, are linearly independent; while for every $i,j\in\{1,\dots,n-1 \}$ it holds $\braket{w_i}{z_j}=0$. Moreover, the following relations holds true:
	\begin{equation}\label{eq:orthogonal_to_1}
		\begin{split}
			\braket{w_i}{\ind_{\cX_1}}=\braket{w_i}{\ind_{\cX_2}}=
			\braket{z_i}{\ind_{\cX_1}}=\braket{z_i}{\ind_{\cX_2}}=0,\quad\forall i\in\{1,\dots,n-1 \}\,,
		\end{split}
	\end{equation}
	We now show that, depending on the mass that the initial distribution $\alpha$ gives to the two blocks $\cX_1$ and $\cX_2$, three different behaviors for the Yaglom limit $\varphi_\star^\alpha$ can take place. For simplicity in what follows we will identify $\cX_1=\{1,\dots,n \}$ and $\cX_2=\{n+1,\dots,2n \}$.
	\begin{proposition}\label{prop:meanfield}
		Depending on the initial distribution $\alpha$ the distribution $\varphi_\star^\alpha$ exhibits the following trichotomy:
		\begin{enumerate}
			\item If $\alpha(\cX_1)>\alpha(\cX_2)$ then $\bra{\varphi_\star^\alpha}=\frac{1}{n}\bra{\ind_{\cX_1}}$.
			\item If $\alpha(\cX_1)<\alpha(\cX_2)$ then $\bra{\varphi_\star^\alpha}=\frac{1}{n}\bra{\ind_{\cX_2}}$.
			\item If $\alpha(\cX_1)=\alpha(\cX_2)$ then 
\begin{equation}\label{eq:limit-phi-equal-size}
	\bra{\varphi_\star^\alpha}=\frac{1}{2n}\[\bra{1}+\frac{1}{\frac1n-2\alpha_{\min}}\(\sum_{i=1}^{n-1}\(\frac1n-2\alpha(i+1) \)\bra{w_i}+\(\frac1n-2\alpha(n+i+1) \)\bra{z_i}\) \]\,.
\end{equation}
		\end{enumerate}
	\end{proposition}
	Before proving \cref{prop:meanfield}, let us briefly comment on its meaning. Items $(1)$ and $(2)$ are quite easy to interpret: if the initial distribution gives a larger mass to one of the two blocks, then, conditionally on not having reached an equilibrium at a very large time, the distribution of the process must be uniformly spread on the block that received a larger initial mass. On the other hand, item $(3)$ represent the ``degenerate'' case in which the two blocks start already with their right equilibrium mass. In this case, the Yaglom limit preserves the proportions between the two blocks; nonetheless, by \cref{eq:halting-alpha} there must be some states out of the support of $\varphi_\star^\alpha$ and \cref{eq:limit-phi-equal-size} states that those coincide with the minimizers of $\alpha$.
	\begin{proof}
		We can write
		\begin{equation}\label{eq:alpha_decomp}
			\bra{\alpha}=\bra{\pi}+a \bra{v}+ \sum_{i=1}^{n-1} b_i\bra{z_i}+c_i \bra{w_i}  \,.
		\end{equation}
	We focus on item $(1)$, the proof of item $(2)$ being completely symmetric. Using the decomposition in \cref{eq:alpha_decomp} and the assumption $\alpha(\cX_1)>\alpha(\cX_2)$ we obtain
	\begin{equation}
		\frac{1}{2}<	\braket{\alpha}{\ind_{\cX_1}}=\frac{1}{2}+a \sqrt{\frac{n}{2}}\,,
	\end{equation}
	from which we deduce that $a>0$, so that, by \cref{eq:varphistar_simplified} we deduce
	\begin{equation}
		\bra{\varphi_\star^\alpha}=\bra{\pi}+ \frac{1}{\max_{y\in\cX}-\frac{v(y)}{\pi(y)}}\bra{v}=\bra{\pi}+\frac{1}{\sqrt{2n}}\bra{v}=\frac{1}{n}\bra{\ind_{\cX_1}}\,.
	\end{equation}
		We now prove $(3)$. Under the assumption $\alpha(\cX_1)=\alpha(\cX_2)$, thanks to \cref{eq:orthogonal_to_1}, we have
		\begin{equation}
			\frac{1}{2}=	\braket{\alpha}{\ind_{\cX_1}}=\frac{1}{2}+a \sqrt{\frac{n}{2}}\,,
		\end{equation}
		from which  it follows that $a=0$. 
Since $\alpha(\cX_1)=\alpha(\cX_2)=\frac{1}{2}$, there exist two probability distributions $\vec p\in\cP(\cX_1)$ and $\vec q\in\cP(\cX_2)$ such that
$$\alpha=\frac{1}{2}\(p_1,\dots,p_n,q_1,\dots,q_n\)\,.$$
Then, called 
$$p_{\min}=\min\{p_1,\dots,p_n\},\qquad q_{\min}=\min\{q_1,\dots,q_n\}\,,$$
and assuming, w.l.o.g., that $q_{\min}\le p_{\min}$ what we have to prove, i.e., \cref{eq:limit-phi-equal-size}, reduces to
\begin{equation}\label{eq:phistar-pq}
	\bra{\varphi_\star^\alpha}=\frac{1}{2n}\(\bra{1}+\frac{1}{\frac1n-q_{n}}\sum_{i=1}^{n-1}\(\frac1n-p_{i+1} \)\bra{w_i}+\(\frac1n-q_{i+1} \)\bra{z_i} \)\,.
\end{equation}
By decomposing $\alpha$ as in \cref{eq:alpha_decomp} with $a=0$, we deduce that
		\begin{align}\label{eq:coefficients}
			c_i=\frac12\(\frac1n-p_{i+1} \),\quad b_i=\frac12\(\frac1n-q_{i+1} \)\,,\qquad i=1,\dots,n-1\,.
		\end{align}
	Therefore, the normalizing constant $\ell^\alpha$ in \cref{eq:ell-alpha} is given by
	\begin{dmath}
		\ell^\alpha=2n \max\left\{c_1,\dots,c_{n-1},b_1,\dots,b_{n-1},-\sum_{i=1}^{n-1}c_i,-\sum_{i=1}^{n-1}b_i \right\}\\
		=2n\max\left\{\frac12\(\frac1n-p_2\),\dots,\frac12\(\frac1n-p_n \),\frac12\(\frac1n-q_2\),\dots,\frac12\(\frac1n-q_n\),\\ -\frac{n-1}{2n}+\frac12(1-p_1),-\frac{n-1}{2n}+\frac12(1-q_1) \right\}\,.
	\end{dmath}
On the one hand, we notice that the last two elements in the latter set are negative, since
$$-\frac{n-1}{2n}+\frac1n(1-p_1)=\frac12\(\frac1n-p_1\)<0,$$
and similarly for the other. On the other hand,
$$\max\left\{\frac12\(\frac1n-p_2\),\dots,\frac12\(\frac1n-p_n \),\frac12\(\frac1n-q_2\),\dots,\frac12\(\frac1n-q_n\)\right\}=\frac12\(\frac1n-q_n\)\ge0\,,$$
therefore
\begin{equation}\label{eq:ell-alpha-blocks}
		\ell^\alpha=n\(\frac1n-q_n\).
\end{equation}
Plugging \cref{eq:coefficients,eq:ell-alpha-blocks} into \cref{eq:phi-star-ell} we get \cref{eq:phistar-pq}.
	\end{proof}

\subsection{Low temperature 1-$d$ random walk in FW regime}\label{sec:1drw}
Let us now consider the following example:
\begin{equation}\label{eq:P-FW}
	P=\left(\begin{matrix}
		1-a&a&0\\
		b&1-b-c&c\\
		0&d&1-d
	\end{matrix} \right),
\end{equation}
where we assume that for some $\beta>0$,
$$a=e^{-\beta\Delta_1}\,,\quad b=e^{-\beta\Delta_2}\,,\quad c=e^{-\beta\Delta_3}\,, \quad d=e^{-\beta\Delta_4}\,.$$
The general case will be treated numerically in \cref{conclusions}.
Here, to make the analysis simpler, we will consider the special case in which the communication height between $3$ and $2$  coincides with that between $1$ and $2$, i.e., $\Delta_1=\Delta_4$ and therefore $a=d$. Moreover, we will assume that
\begin{equation}\label{eq:ch}
	\Delta_2>\Delta_1=\Delta_4>\Delta_3\,,
\end{equation}
so that, for $\beta\gg 1$
$$0<b\ll a=d\ll c\ll 1 \,. $$
\begin{figure}[h]
	\includegraphics[height=6cm]{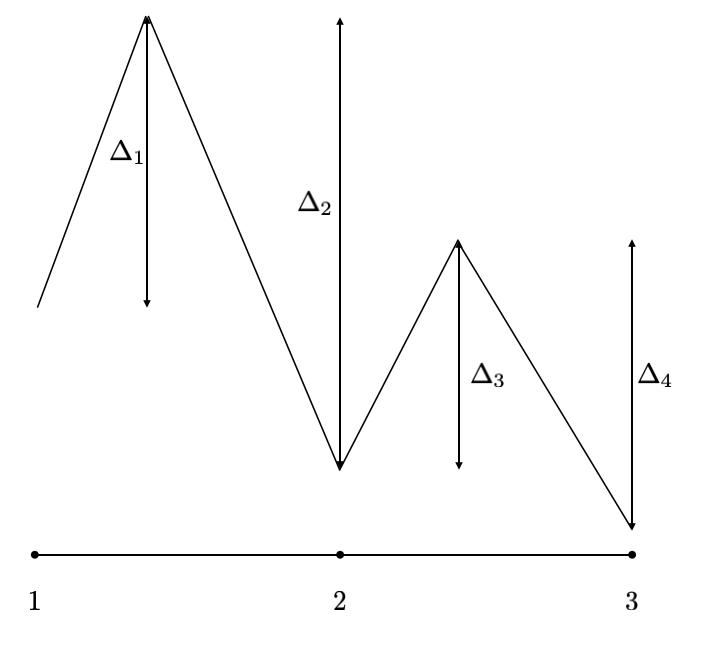}
	\caption{The energy landscape associated with the transition matrix in \cref{eq:P-FW} under the condition in \cref{eq:ch}.}
\end{figure}
It is immediate to check that
$$\pi=\frac{1}{Z}\(\frac{b}{c},\frac{a}{c},1\)\,, $$
where
$$ Z=\frac{a+b+c}{c}\,.$$
Moreover, we have also the following eigenpairs
$$\lambda_2=1-a,\qquad \nu_2=(1,0,-1)\,,$$
and
$$\lambda_3=1-a-b-c,\qquad \nu_3=\(-\frac{b}{b+c},1,-\frac{c}{b+c}\)\,.$$
Therefore, for any probability distribution $\alpha$ we have
$$\bra{\alpha}=\bra{\pi}+a_2\bra{\nu_2}+a_3\bra{\nu_3}\,,$$
where the coefficients $a_2$ and $a_3$ are uniquely determined by $\alpha$ through the relation
\begin{align*}
	a_2&=\alpha(1)-\frac{b}{cZ}+\frac{b}{b+c}\(\alpha(2)-\frac{a}{cZ} \)\\
	a_3&=\alpha(2)-\frac{a}{cZ}\,.
\end{align*}
Therefore, for every $\alpha$ there are only four available options for the associated Yaglom limit, depending on the signs of $a_2$ and $a_3$:
\begin{figure}[h]
	\includegraphics[height=7cm]{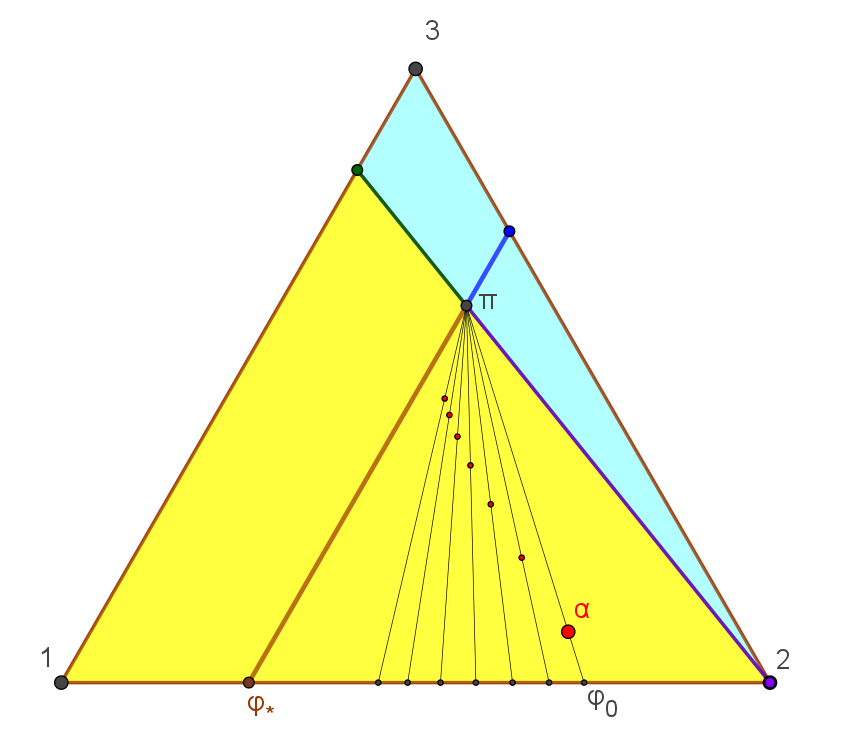}
	\includegraphics[height=7cm]{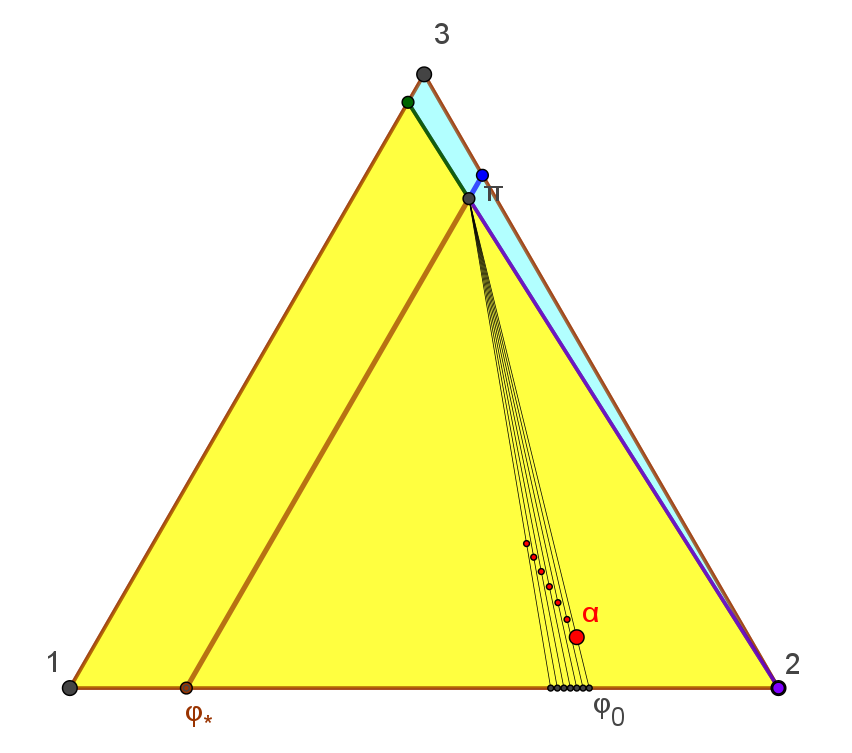}
	\caption{	A numerical plot of the basins of attraction for the model in \cref{sec:1drw}, in the case $a=d$. The yellow region represents the basin of attraction of the Yaglom limit in the bottom edge of the simplex (the brown thick point), while the the cyan region is the basin of attraction of the blue thick point on the right edge of the simplex. As explained, there are other two Yaglom limits (the purple and the green thick points), but their basin of attraction is lower dimensional (the purple and the green segments, respectively). We consider the evolution of the process started at $\alpha$ (the red thick point in the yellow region): the successive red points represent the evolution of the measure $\mu_t^\alpha$ for $t=0,\dots,6$. The projection of these points on the edges of the simplex (the black points), represent the corresponding conditional measures $\varphi^\alpha_t$ for $t=0,\dots, 6$. As this picture suggests, the distribution $\mu_t^\alpha$ will approach the stationary distribution $\pi$ by the the direction of the brown segment which links $\pi$ to the Yaglom limit $\varphi^\star_\alpha$. The difference between the two plots is only due to a change in the inverse temperature $\beta$, which is higher in the second plot.
	}
\end{figure}

\begin{enumerate}
	\item $a_2>0$: In this case we get 
	$$\ell^{\alpha}=\max_{y} -a_2\frac{\nu_2(y)}{\pi(y)}=a_2 Z\,,$$ and
	\begin{equation}\label{case1}
	\bra{\varphi^\alpha_{\star}}=\bra{\pi}+\frac{a_2}{\ell^{\alpha}}\bra{\nu_2}=\(\frac{b+c}{a+b+c},\frac{a}{a+b+c},0 \)\,.
	\end{equation}
	\item $a_2<0$: Then	\begin{align*}
		\ell^{\alpha}&=\max_{y} -a_2\frac{\nu_2(y)}{\pi(y)}=-a_2\frac{Z}{bc} \,,
	\end{align*}
	and
	\begin{align}\label{case2}	\bra{\varphi^\alpha_{\star}}&=\bra{\pi}+\frac{a_2}{\ell^{\alpha}}\bra{\nu_2}=\(0,\frac{a}{a+b+c},\frac{b+c}{a+b+c}\)\,.
	\end{align}
	\item $a_2=0$ and $a_3>0$: We get, 
	\begin{align*}
		\ell^{\alpha}&=\max_{y} -a_3\frac{\nu_3(y)}{\pi(y)}=\(1-\frac{a}{cZ}\)\frac{cZ}{b+c} \,,
	\end{align*}
and
\begin{align}\label{case3}	\bra{\varphi^\alpha_{\star}}&=\bra{\pi}+\frac{a_3}{\ell^{\alpha}}\bra{\nu_3}=(0,1,0)\,.
\end{align}
	\item $a_2=0$ and $a_3<0$: In this case	\begin{align*}
		\ell^{\alpha}&=\max_{y} -a_3\frac{\nu_3(y)}{\pi(y)}=-a_3\frac{Z}{ac} \,,
	\end{align*}
	and
	\begin{align}\label{case4}	\bra{\varphi^\alpha_{\star}}&=\bra{\pi}+\frac{a_3}{\ell^{\alpha}}\bra{\nu_3}=\(\frac{b}{b+c},0,\frac{c}{b+c}\)\,.
	\end{align}
\end{enumerate}
Among the four possible Yaglom limits there is a prominent one, that in \cref{case1}, the basin of attraction of which invades  essentially the whole  simplex in the low temperature limit (in which the stationary distribution becomes concentrated on state $3$). In the latter case is then natural to compare the distribution in \cref{case1} with the classical notion of quasi-stationary distribution, that is, when the state $3$ is made absorbing.
At this scope, we consider the restricted matrix
\begin{eqnarray}
	[P]=\left(\begin{matrix}
		1-a&a\\
		b&1-b-c
	\end{matrix} \right)\,,
\end{eqnarray}
and by simple spectral analysis it is easy to check that the process admits a unique quasi-stationary distribution
\begin{equation}\label{absorbing}
	\begin{split}
	\phi^*&=\left(\frac{\sqrt{2 c (b-a)+(a+b)^2+c^2}-a-b+c}{2 c},\frac{2 a}{\sqrt{2 c
			(b-a)+(a+b)^2+c^2}+a+b+c}\right)\\
		&\sim \(\frac{c}{a+c},\frac{a}{a+c} \)\,,
		\end{split}
\end{equation}
where the asymptotic notation refers to the limit $\beta\to\infty$. The 
associated eigenvalue is
$$\lambda^*=\frac{1}{2} \left(\sqrt{2 c (b-a)+(a+b)^2+c^2}-a-b-c+2\right)\sim1-a \,.$$
We deduce that in the low temperature limit our notion of quasi-stationary distribution is well-approximated by the classical notion, since $\lambda^*\sim\lambda_2$ and the two distributions in \cref{case1} and \cref{absorbing} are asymptotically equivalent in a multiplicative sense. We will comment further on the analogy between absorbing and non-absorbing case in \cref{conclusions}.

\section{Counterexamples in the non-reversible setting}
In this section we present two examples of non-reversible ergodic Markov chains in which there exist some initial distributions $\alpha$ such that the Yaglom limit $\varphi_\star^\alpha$ does not exists. We start in \cref{suse:cex1} by providing an example of a three-state chain in which, choosing $\alpha$ has the Dirac distribution on a state, the conditional distribution $(\varphi_t^\alpha)_{t \ge 0}$ converges toward a period cycle of length $4$. Subsequently, in \cref{suse:cex2} we exhibit a four-state chain having a genuine metastable behavior, in which the sequence of conditional distributions $(\varphi_t^\alpha)_{t \ge 0}$ converges toward a periodic cycle that can be easily interpreted.

\subsection{A three-state non-reversible chain}\label{suse:cex1}
Let $p\in(0,1]$ and consider the transition matrix
\begin{equation}
	P=\frac13\left(\begin{matrix}
		1&1+p&1-p\\
		1-p&1&1+p\\
		1+p&1-p&1
	\end{matrix} \right)=\ket{1}\bra{\pi}+\frac{p}{\sqrt{3}}A\,,
\end{equation}
where $\pi$ is the uniform distribution (and, by the bistochasticity of $P$, is the unique stationary distribution of $P$) and
\begin{equation}
	A=\frac{1}{\sqrt{3}}\left(
\begin{matrix}
	0&1&-1\\
	-1&0&1\\
	1&-1&0
\end{matrix}	
	 \right)\,.
\end{equation}
One can easily compute ${\rm spec}(P)=(1,\pm\tfrac{p\, \mathbf{i}}{\sqrt{3}})$. Note that for all $t>0$
\begin{equation}
	A^t=\begin{cases}
		A&t\equiv 1\ {\rm mod}(4)\,,\\
		B&t\equiv 2\ {\rm mod}(4)\,,\\
		-A&t\equiv 3\ {\rm mod}(4)\,,\\
		-B&t\equiv 0\ {\rm mod}(4)\,,
	\end{cases}
\end{equation}
where
\begin{equation}
	B=\frac{1}{3}\left(
	\begin{matrix}
		-2&1&1\\
		1&-2&1\\
		1&1&-2
	\end{matrix}	
	\right)\,.
\end{equation}
Then
\begin{equation}
	P^t=\ket{1}\bra{\pi}+\(\frac{p}{\sqrt{3}}A\)^t+\ket{1}\bra{\pi}\sum_{s=1}^{t-1}\binom{t}{s}\(\frac{p}{\sqrt{3}}A\)^s\,.
\end{equation}
Therefore, for any $y\in\{1,2,3\}$, and $t>0$
\begin{align*}
	P^t(y,\cdot)-\pi(\cdot)&=\(\frac{p}{\sqrt{3}}A\)^t(y,\cdot)-\sum_x\pi(x) \sum_{s=1}^{t-1}\binom{t}{s}\(\frac{p}{\sqrt{3}}A\)^s(x,\cdot)\\
	&=\(\frac{p}{\sqrt{3}}A\)^t(y,\cdot)-\frac13 \sum_{s=1}^{t-1}\binom{t}{s}\sum_x\(\frac{p}{\sqrt{3}}A\)^s(x,\cdot)\\
	&=\(\frac{p}{\sqrt{3}}A\)^t(y,\cdot)\,,
\end{align*}
where in the last line we used that $\sum_xA^k(x,y)=0$ for all $y\in\{1,2,3\}$. Let $\alpha=\delta_1$. and notice that
\begin{equation}
	s_t^\alpha=\max_x 1-\frac{\mu_t^\alpha(x)}{\pi(x)}=1-3\min_x \mu_t^\alpha(x)\,,\qquad t\ge 0\,.
\end{equation}
Therefore
\begin{equation}
	s_t^\alpha=\begin{cases}
		1&t=0\,,\\
		1-3\(\frac{1}{3}-\frac{1}{\sqrt{3}}\(\frac{p}{\sqrt{3}}\)^t\)=\sqrt{3}\(\frac{p}{\sqrt{3}}\)^t&t\equiv 1\ {\rm mod}(4)\,,\\
	1-3\(\frac{1}{3}-\frac{2}{3}\(\frac{p}{\sqrt{3}}\)^t\)=2\(\frac{p}{\sqrt{3}}\)^t&t\equiv 2\ {\rm mod}(4)\,,\\
		1-3\(\frac{1}{3}-\frac{1}{\sqrt{3}}\(\frac{p}{\sqrt{3}}\)^t\)=\sqrt{3}\(\frac{p}{\sqrt{3}}\)^t&t\equiv 3\ {\rm mod}(4)\,,\\
	1-3\(\frac{1}{3}-\frac{1}{3}\(\frac{p}{\sqrt{3}}\)^t\)=\(\frac{p}{\sqrt{3}}\)^t&t\equiv 0\ {\rm mod}(4)\,.
	\end{cases}
\end{equation}
Recall from \cref{eq:limiting-ratio} that the existence of the Yaglom limit is equivalent to the existence of the limiting vector in \cref{eq:limiting-ratio}. In this example we have for all $t>0$
\begin{equation}
	\frac{1}{s_t^\alpha}\left(\bra{\mu_t^\alpha}-\bra{\pi}\right)=\begin{cases}
		\frac{1}{\sqrt{3}}A(1,\cdot)=(0,\tfrac13,-\tfrac13)&t\equiv 1\ {\rm mod}(4)\,,\\
		\frac{1}{2}B(1,\cdot)=(-\tfrac13,\tfrac16,\tfrac16)&t\equiv 2\ {\rm mod}(4)\,,\\
		-\frac{1}{\sqrt{3}}A(1,\cdot)=(0,-\tfrac13,\tfrac13)&t\equiv 3\ {\rm mod}(4)\,,\\
-B(1,\cdot)=(\tfrac23,-\tfrac13,-\tfrac13)&t\equiv 4\ {\rm mod}(4)\,,
		\end{cases}
\end{equation}
thus, the limit does not exists. Notice that we have four  possible limits along subsequences for $(\varphi^\alpha_t)_{t\ge 0}$, being
\begin{align*}
	(\tfrac13,\tfrac23,0)\,,\qquad  (0,\tfrac12,\tfrac12)\,,\qquad	(\tfrac13,0,\tfrac23)\,,\qquad  (1,0,0)\,.
\end{align*}
\subsection{A four-state non-reversible metastable system}\label{suse:cex2}
Fix $\varepsilon\in(0,1)$ and consider the transition matrix
\begin{equation}
	P=\left(\begin{matrix}
		0&1-\varepsilon&0&\varepsilon\\
		0&0&1-\varepsilon&\varepsilon\\
		1-\varepsilon&0&0&\varepsilon\\
		\frac{1-\varepsilon}{3}&\frac{1-\varepsilon}{3}&\frac{1-\varepsilon}{3}&\varepsilon
	\end{matrix} \right)\,.
\end{equation}
Notice that $\pi$, i.e., the unique stationary distribution of $P$, coincides with its last row. Notice that ${\rm spec}(P)=(1,0,-\tfrac12(1-\varepsilon)(1\pm\sqrt{3}\ \mathbf{i}))$.
Let $\alpha=\delta_1$.
The behavior of the system is easy to describe: the process will evolve periodically along the cycle $1\to2\to3\to 1$ up to a geometric time of rate $\varepsilon$, when it hits the state $4$, and at the next step it is a equilibrium. When $\varepsilon\ll 1$, this is a typical metastable situation, but the metastable state is described by a periodic orbit on the simplex, rather than by a single measure. As the physical intuition suggests, we now show that the sequence of conditional distributions $(\varphi_t^\alpha)_{t\ge 0}$ does not converge, but it is constant along subsequences depending on the equivalence class of $t$ modulo $3$.
As first step, it is immediate to check that $s_0^\alpha=1$, henceforth $\varphi_0^\alpha=\alpha$. Similary, being $P(1,\cdot)$ non-fully supported, $s_1^\alpha=1$ and therefore $\varphi_1=P(1,\cdot)$. Now notice that
\begin{eqnarray*}
	P^2(1,\cdot)=\left(\tfrac13\varepsilon(1-\varepsilon)\,,\tfrac13\varepsilon(1-\varepsilon)\,, \tfrac13(3-2\varepsilon)(1-\varepsilon)\,, \varepsilon \right)\,,
\end{eqnarray*}
so that $s_2^\alpha=1-\varepsilon$. By using \cref{lemma:recursion-phi} we deduce $\varphi_2^\alpha=P(2,\cdot)$. In general, one can check that for all $t\ge 2$
\begin{equation}
	\min_x \frac{P^t(1,\cdot)}{\pi(x)}= \frac{(t-1)\varepsilon+\sum_{i=2}^t(-1)^{i+1}\binom{t}{i}\varepsilon^i}{1-\varepsilon}\,,
\end{equation}
from which follows that
\begin{equation}
	s_t=(1-\varepsilon)^{t-1}\,,\qquad t\ge 2\,,
\end{equation}
and therefore
\begin{equation}\label{eq:ratio-sepa}
	\frac{s_{t+1}}{s_t}=1-\varepsilon\,,\qquad t\ge 1\,.
\end{equation}
In conclusion, we obtain
\begin{equation}
	\varphi_t^\alpha=\begin{cases}
		(1,0,0,0)&t=0\,,\\
		(0,1-\varepsilon,0,\varepsilon)&t\equiv1\ {\rm mod}(3)\,,\\
		(0,0,1-\varepsilon,\varepsilon)&t\equiv2\ {\rm mod}(3)\,,\\
		(1-\varepsilon,0,0,\varepsilon)&t\equiv0\ {\rm mod}(3)\,.
	\end{cases}
\end{equation}
Moreover, by \cref{lemma:recursion-phi,eq:ratio-sepa} we also have
\begin{equation}
	\bra{\varphi_{t}^\alpha} P=(1-\varepsilon)\bra{\varphi_{t+1}^\alpha} +\varepsilon \bra{\pi}\,,\qquad t\ge 1\,.
\end{equation}

\section{Conclusions and future research}\label{conclusions}
In this work we proposed a new notion of quasi-stationary distribution for non-absorbing reversible Markov chains.
In particular we characterized such distributions both by a suitably defined
Yaglom limit, and by relating it to spectral quantities.  The geometric interpretation of our results represents a new point of view on the evolution of the conditioned process and its dependence on the starting distribution $\alpha$.

As mentioned in the introduction, our main motivation is a general characterization of metastability,
deriving asymptotic estimates in the framework of an exact representation formula, in the spirit of the result presented in \cite[Theorem 3]{MS19} for the absorbing case.
Nevertheless, such an extension to non-absorbing chains is far from trivial and represents our main
open problem at this step. 
More precisely, in \cite{MS19}  a notion of
\emph{conditional strong quasi-stationary time} was introduced to describe the local relaxation time and to provide
the above-mentioned exact representation formula, from which rich asymptotic results can be obtained, especially on the
dependence of the asymptotic exponential law on the initial distribution. 
\begin{figure}[h]
		\includegraphics[height=6.5cm]{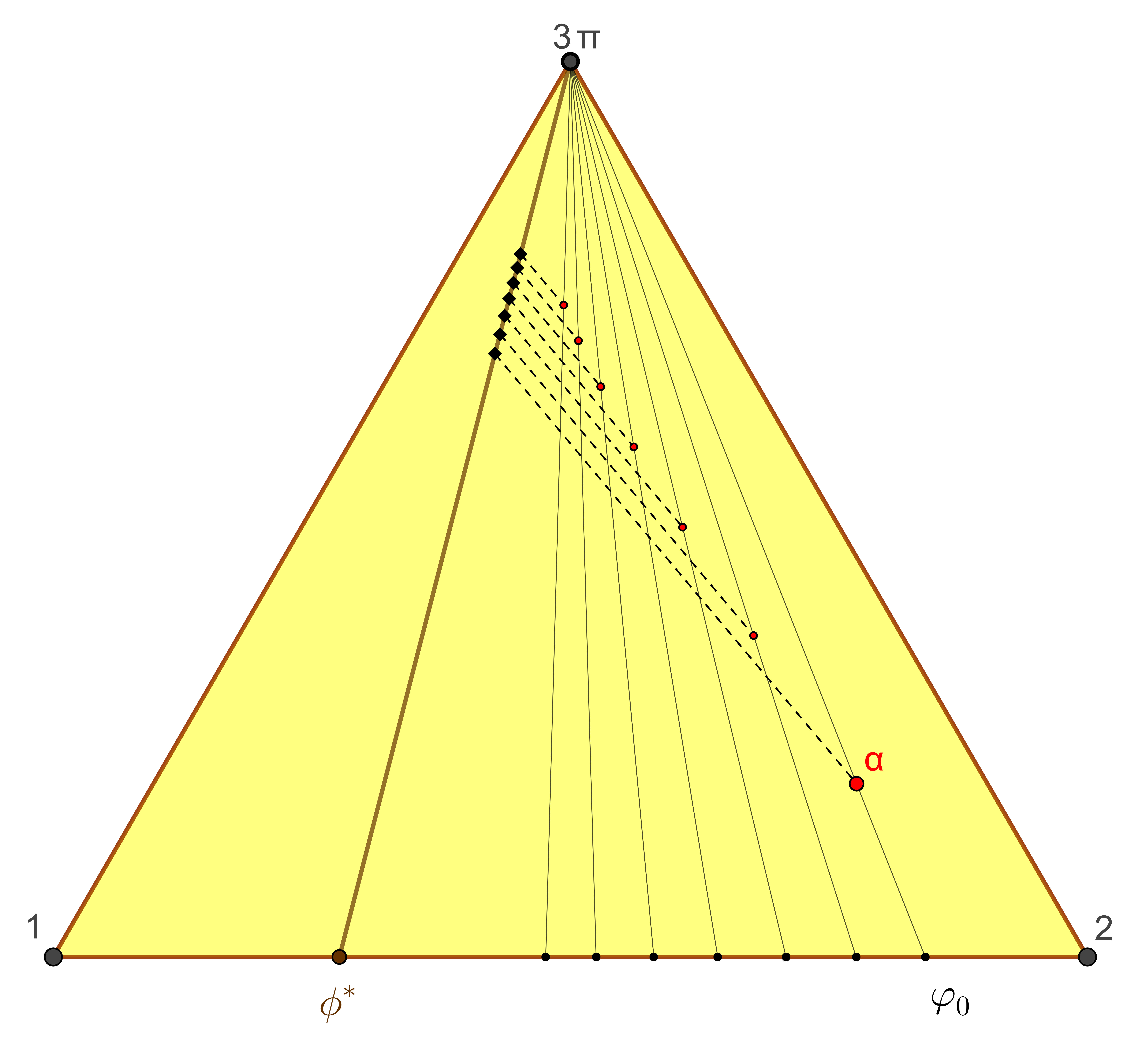}\qquad
	\includegraphics[height=6.5cm]{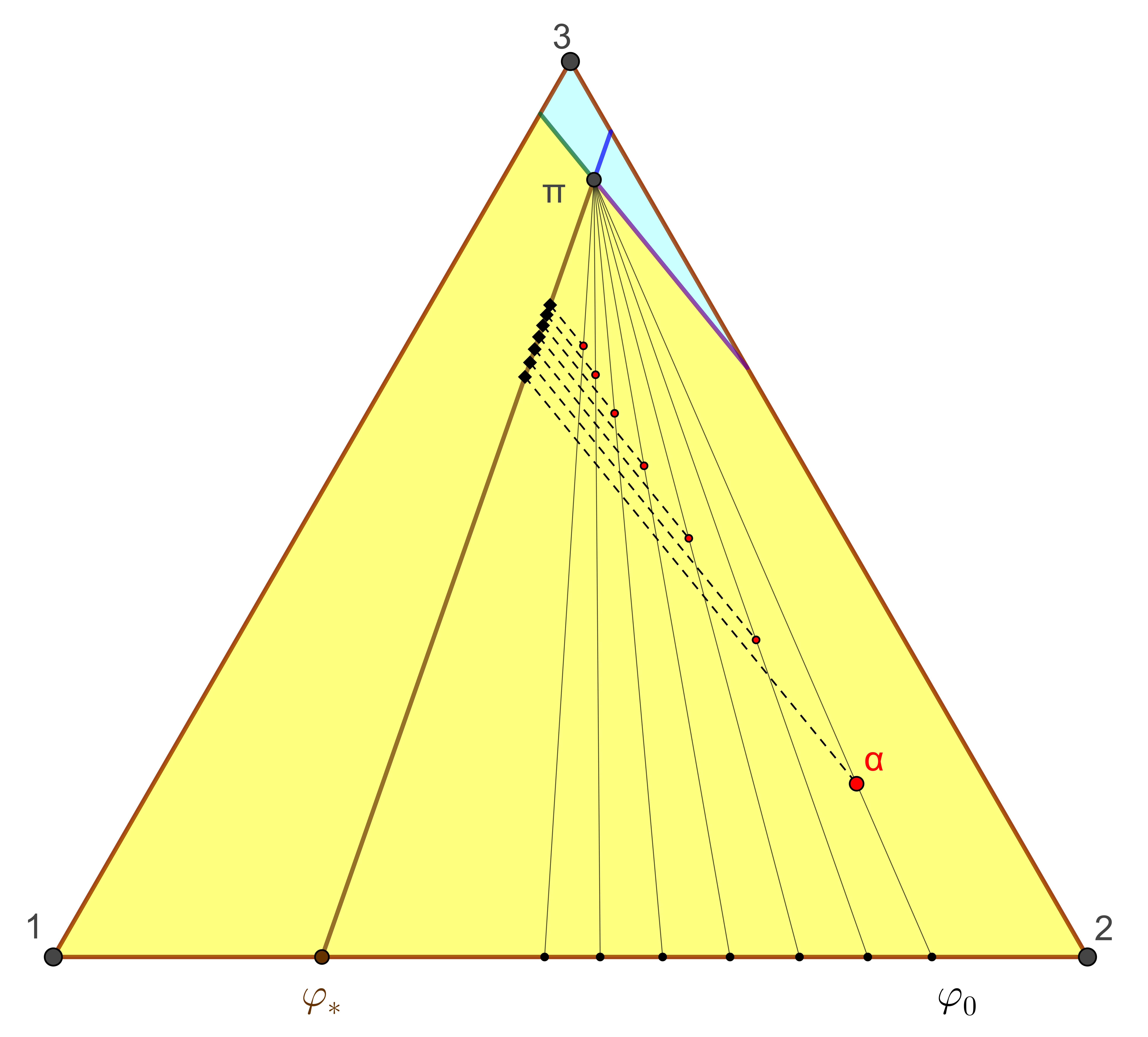}
	\caption{	{(A) A numerical simulation of the example presented in \cref{sec:1drw} (in the general setup $a\neq d$) in the absorbing case $d=0$.
	The red dots represent the evolution of the process $\mu^\alpha_t$, while the black squares represent the evolution of
	the process $\mu^{\bar\alpha}_t$.
	(B) A numerical simulation of the same model as in (A), but in the non-absorbing case, $d\not=0$, with the same choice  as in (A) for the other parameters.
		}
	}\label{fig:abs}
\end{figure}

In order to explain the difficulties of the non-absorbing case, let us present the ideas of \cite{MS19} using the notation and the geometric language developed in the current paper.  In \cref{fig:abs}(A) we represent an absorbing process. The key idea to obtain the estimates in \cite{MS19} is to quantify the convergence of the sequence $(\mu^\alpha_t)_{t\ge 0}$ (represented by the red dots in \cref{fig:abs}) to the sequence $(\mu_t^{\bar{\alpha}})_{t\ge 0}$ (represented by the black squares in \cref{fig:abs}) defined as follows
\begin{equation}\label{eq:beta-1}
	\begin{split}
		\bra{\mu^{\bar{\alpha}}_t}\coloneqq\bra{\bar{\alpha}}P^t&=\braket{\mu^\alpha_t}{f_2} \bra{\phi^*}+(1- \braket{\mu^\alpha_t}{f_2})\bra{\pi}\\
		&=\lambda_2^t\braket{\alpha}{f_2} \bra{\phi^*}+(1- \lambda_2^t\braket{\alpha}{f_2})\bra{\pi}\,,
		\end{split}
\end{equation}
where 
\begin{equation}
	\bra{\bar{\alpha}}\coloneqq \braket{\alpha}{f_2} \bra{\phi^*}+(1- \braket{\alpha}{f_2})\bra{\pi}\,,
	\end{equation}
	coincides with the projection (parallel to $v_3$) of $\alpha$ on the line passing through $\pi$ and $\phi^*$. It is crucial to observe that the sequence $(\mu^{\bar{\alpha}}_t)_{t\ge 0}$ is all contained in the latter line and its separation from $\pi$ is trivial thanks to the second representation in \cref{eq:beta-1}. To control the convergence of $\mu_t^\alpha$ to $\mu_t^{\bar{\alpha}}$ the authors in \cite{MS19} introduce a local notion of separation distance, which is then controlled by the help of a ``local chain'' (on the set of non-absorbing states, $\cM$) related to the Doob transform
	of the  transition matrix restricted to ${\cM}$. The convergence to equilibrium of this local chain,  in separation distance, is then used to control the distribution of the \emph{conditional strong quasi-stationary time}, see \cite[Theorem 1]{MS19}. 

In figure  \cref{fig:abs}(B) the evolution of the non absorbing process is shown.
Despite the similarity of the behaviour of the processes, the extension of the analysis to the non-absorbing case 
requires a new notion of ``locality'' which is the main open problem at this stage. More precisely it is not
clear how to show the existence of an analogue \emph{conditional strong quasi-stationary time}, how to define a local separation distance
and how to generalize the definition of the local chain.

Another  interesting open problem is the generalization of the results to non-reversible chain.
As discussed in the four-state non reversible example the conditional distributions $\phi_t^\alpha$ may not to
converge to a Yaglom limit, but to a cyclic orbit in the simplex of probability measures. This example suggests
 very interesting new scenarios in the study of the behaviour of the conditioned distributions $\phi^\alpha_t$,
 connected with non trivial behaviour of dynamical systems.

 \end{document}